\documentclass[12pt,a4paper]{article}
\usepackage{indentfirst}
\usepackage{amsmath,amssymb,amsfonts,amsthm,graphics}
\usepackage{makeidx}
\usepackage{amsmath,amssymb,amsfonts,amsthm,graphics,graphicx}
\usepackage{makeidx}
\usepackage{color}
\usepackage{ifpdf}
\usepackage{subfigure}

\newtheorem{thm}{Theorem}

\newtheorem{lem}{Lemma}
\newtheorem{pro}{Proposition}

\theoremstyle{definition}

\def\-{\mbox{--}}

\def\pf{\noindent {\it Proof.} }
\begin{document}

\title{\Large\bf Proper vertex connection and graph operations}
\author{\small Yingying Zhang, Xiaoyu Zhu\\
\small Center for Combinatorics and LPMC\\
\small Nankai University, Tianjin 300071, China\\
\small E-mail: zyydlwyx@163.com; zhuxy@mail.nankai.edu.cn}
\date{}
\maketitle
\begin{abstract}

A path in a vertex-colored graph is a {\it vertex-proper path} if
any two internal adjacent vertices differ in color. A vertex-colored
graph is {\it proper vertex $k$-connected} if any two vertices of
the graph are connected by $k$ disjoint vertex-proper paths of the
graph. For a $k$-connected graph $G$, the {\it proper vertex
$k$-connection number} of $G$, denoted by $pvc_{k}(G)$, is defined
as the smallest number of colors required to make $G$ proper vertex
$k$-connected. A vertex-colored graph is {\it strong proper
vertex-connected}, if for any two vertices $u,v$ of the graph, there
exists a vertex-proper $u$-$v$ geodesic. For a connected graph $G$,
the {\it strong proper vertex-connection number} of $G$, denoted by
$spvc(G)$, is the smallest number of colors required to make $G$
strong proper vertex-connected. In this paper, we study the proper vertex $k$-connection number and the strong proper vertex-connection number on the join of two graphs, the Cartesian, lexicographic, strong and direct product, and present exact values or upper bounds for these operations of graphs.


\end{abstract}

\section{Introduction}

In this paper, all graphs considered are simple, finite and undirected. We refer to the
book \cite{B} for undefined notation and terminology in graph theory. For simplicity,
a set of internally vertex-disjoint paths will be called {\it disjoint}. A path in
an edge-colored graph is a {\it rainbow path} if its edges have different colors. An
edge-colored graph is {\it rainbow $k$-connected} if any two vertices of the graph
are connected by $k$ disjoint rainbow paths of the graph. For a $k$-connected graph
$G$, the {\it rainbow $k$-connection number} of $G$, denoted by $rc_{k}(G)$, is defined
as the smallest number of colors required to make $G$ rainbow $k$-connected. This concept came
after the September 11, 2001 terrorist attacks, from which many weaknesses of the secure
information transfer of USA had been discovered. Among all these weaknesses, the most nonnegligible one may be the problem of the secure communication between departments of government. There is an
information transfer path (may pass by some other departments) between every pair of departments and each step of this path needs a password. In order to protect the safety system from the invasion of terrorists, all passwords in the path must be different. As the data size can be quite huge, one natural question arose that what is the smallest number of passwords allowed to ensure one or more secure paths between every pair of departments. This concept
was first introduced by Chartrand et al. in \cite{CJM, CJMZ}. Since then, a lot of
results on the rainbow connection have been obtained; see \cite{GMP,LSS, LSu}.

As a natural counterpart of the concept of rainbow $k$-connection, the concept of
rainbow vertex $k$-connection was first introduced by Krivelevich and Yuster
in \cite{KY} for $k=1$, and then by Liu et al. in \cite{LMS} for general $k$.
A path in a vertex-colored graph is a {\it vertex-rainbow path} if its internal
vertices have different colors. A vertex-colored graph is {\it rainbow vertex $k$-connected}
if any two vertices of the graph are connected by $k$ disjoint vertex-rainbow paths of the
graph. For a $k$-connected graph $G$, the {\it rainbow vertex $k$-connection number} of
$G$, denoted by $rvc_k(G)$, is defined as the smallest number of colors required to make
$G$ rainbow vertex $k$-connected. There are many results on this topic, we refer to \cite{CLS1,CLS2,LS,MYWY}.

Concerning about the geodesics instead of the paths, Li et al. \cite{LMS1} introduced the concept of strong rainbow
vertex-connection. A vertex-colored graph is {\it strong rainbow
vertex-connected}, if for any two vertices $u, v$ of the graph,
there exists a vertex-rainbow $u$-$v$ geodesic, i.e., a $u$-$v$ path
of length $d(u,v)$. For a connected graph $G$, the {\it strong
rainbow vertex-connection number} of $G$, denoted by $srvc(G)$, is
the smallest number of colors required to make $G$ strong rainbow
vertex-connected.

In 2011, Borozan et al. \cite{BFG} introduced the concept of proper
$k$-connection of graphs. A path in an edge-colored graph is a {\it
proper path} if any two adjacent edges differ in color. An
edge-colored graph is {\it proper $k$-connected} if any two vertices
of the graph are connected by $k$ disjoint proper paths of the
graph. For a $k$-connected graph $G$, the {\it proper $k$-connection
number} of $G$, denoted by $pc_{k}(G)$, is defined as the smallest
number of colors required to make $G$ proper $k$-connected. This concept is also
based on the situation we introduce above if we bring down a little bit our demand.
That is, we only need to set them different for adjacent passwords instead of any pair of passwords in this path. Note
that $$1\leq pc_k(G)\leq \min\{\chi'(G), rc_k(G)\},\ \ \ \  \ \ \ \
\ \ \ \ \ \  \ \  (1)$$ where $\chi'(G)$ denotes the edge-chromatic
number. Recently, the case for $k=1$ has been studied by Andrews et
al. \cite{ALL}, Laforge et al. \cite{LLZ} and Mao et al. \cite{MYWY1}.

Inspired by the concepts above, Jiang et al. \cite{JLZZ} introduced the concepts of
proper vertex $k$-connection and strong proper vertex-connection. A path in a vertex-colored graph is a
{\it vertex-proper path} if any two internal adjacent vertices
differ in color. A vertex-colored graph is {\it proper vertex
$k$-connected} if any two vertices of the graph are connected by $k$
disjoint vertex-proper paths of the graph. For a $k$-connected graph
$G$, the {\it proper vertex $k$-connection number} of $G$, denoted
by $pvc_{k}(G)$, is defined as the smallest number of colors
required to make $G$ proper vertex $k$-connected. Let $\kappa(G)=$ max\{$k:G$ is $k$-connected\}
denote the vertex-connectivity of $G$. Note that $pvc_{k}(G)$ is
well defined if and only if $1\leq k\leq\kappa(G)$. We write $pvc(G)$
for $pvc_{1}(G)$, and similarly, $rc(G), rvc(G)$ and $pc(G)$ for
$rc_1(G), rvc_1(G)$ and $pc_1(G)$ respectively. For a complete graph $G$, set $pvc(G)=0$. Moreover, we have $pvc(G)\geq 1$ if $G$ is a noncomplete graph. For $k\geq 2$, by definition we
have $pvc_k(G)\geq 1$ if $G$ is a $k$-connected graph.
It is easy to see that $$0\leq pvc_k(G)\leq
\min\{\chi(G), rvc_k(G)\},\ \ \ \  \ \ \ \ \ \ \ \ \ \  \ \  (2)$$
where $\chi(G)$ denotes the chromatic number of $G$. A vertex-colored graph is {\it strong proper
vertex-connected}, if for any two vertices $u, v$ of the graph,
there exists a vertex-proper $u$-$v$ geodesic. For a connected graph
$G$, the {\it strong proper vertex-connection number} of $G$,
denoted by $spvc(G)$, is the smallest number of colors required to
make $G$ strong proper vertex-connected. Especially, set $spvc(G)=0$ for a complete graph $G$. Furthermore, we have $spvc(G)\geq 1$ if $G$ is not complete. Note that if $G$ is a
nontrivial connected graph, then $$0\leq pvc(G)\leq spvc(G)\leq
\min\{\chi(G), srvc(G)\}.\ \ \ \  \ \ \ \ \ \ \ \ \ \  \ \  (3)$$

We recall some fundamental results on $pvc_k(G)$ and $spvc(G)$ which can be found in \cite{JLZZ}.

\begin{pro}\label{pro1}\cite{JLZZ} Let $G$ be a nontrivial connected graph. Then

$(a)$  \ $pvc(G)=0$ if and only if $G$ is a complete graph;

$(b)$ \ $pvc(G)=1$ if and only if $diam(G)=2$.
\end{pro}

For the case that $diam(G)\geq 3$, we have the following theorem.

\begin{thm}\label{thm1}\cite{JLZZ} Let $G$ be a nontrivial connected graph. Then, $pvc(G)=2$ if
and only if $diam(G)\geq 3$.
\end{thm}

Moreover, Jiang et al. \cite{JLZZ} determined the value of $pvc_k(G)$ when $G$ is a complete graph and a complete bipartite graph.

\begin{lem}\label{lemJ1}\cite{JLZZ} $(1)$ $pvc_2(K_n)=pvc_3(K_n)=...=pvc_{n-1}(K_n)=1$.

$(2)$ $pvc_k(K_{n_1,n_2})=2$ for $2\leq k\leq n_1\leq{n_2}$.

\end{lem}

The following results on $spvc(G)$ are immediate from its definition.

\begin{pro}\label{pro2}\cite{JLZZ} Let $G$ be a nontrivial connected graph of order $n$. Then

$(a)$ \ $spvc(G)=0$ if and only if $G$ is a complete graph;

$(b)$ \ $spvc(G)=1$ if and only if $diam(G)=2$.

\end{pro}

The standard products (Cartesian, direct, strong, and lexicographic) draw a constant attention of graph research community, see some papers \cite{ACKP,BS,GV,KS,NS,P,S,Z}. In this paper we consider the join and the four standard products with respect to the (strong) proper vertex-connection number. Each of them will be treated in one of the forthcoming sections. In the join part, we determine the values of the
proper vertex $k$-connection number and the strong proper vertex-connection
number for the join of two graphs. Besides, for the Cartesian, the lexicographic
and the strong products, we also study the two parameters, giving exact values for
most of our results and upper bounds for the others. In the final section, we determine the values of the proper vertex connection number for the direct product, and study the proper vertex $2$-connection number and the strong proper vertex-connection
number for the direct product with one of its factors being the complete graph. For all graph products, only $k=1,2$ are considered in this paper.

\section{The join}

The {\it join} $G\vee H$ of two graphs $G$ and $H$ has vertex set $V(G)\cup V(H)$ and its edge set consists of $E(G)\cup E(H)$ and the set $\{uv:u\in V(G)$ and $v\in V(H)\}$.

\begin{thm}\label{thmJ1} $(1)$ If $G$ and $H$ are graphs such that $G\vee H$ is not complete, then $pvc(G\vee H)=spvc(G\vee H)=1$.

$(2)$ Let $G,H$ be two graphs and $2\leq k\leq \min\{|G|,|H|\}$. If the sum of the minimum degrees of $G$ and $H$ is less than $k-1$, then $pvc_k(G\vee H)=2$; otherwise, we have $pvc_k(G\vee H)=1$.
\end{thm}

\pf (1) By the definition of join, we have $diam(G\vee H)=2$ since $G\vee H$ is not complete. From Propositions \ref{pro1} and \ref{pro2}, it follows that $pvc(G\vee H)=spvc(G\vee H)=1$.

(2) Let $u$ and $v$ be two vertices with the minimum degree in $G$ and $H$, respectively. If the sum of the minimum degrees of $G$ and $H$ is less than $k-1$, then there must exist a path of length at least 3 among the $k$ desired paths from $u$ to $v$ in $G\vee H$. Thus $pvc_k(G\vee H)\geq2$.
Clearly, $G\vee H$ has a spanning complete bipartite subgraph. By Lemma \ref{lemJ1}(2), we have $pvc_k(G\vee H)\leq2$ and so $pvc_k(G\vee H)=2$. For the other cases, we can always find  $k$ desired paths of length at most 2 between any two vertices of $G\vee H$. Thus $pvc_k(G\vee H)=1$.\qed

\section{The Cartesian product}

The {\it Cartesian product} $G\square H$ of two graphs $G$ and $H$ is the graph with vertex set $V(G)\times V(H)$, in which two vertices $(g,h)$ and $(g',h')$ are adjacent if and only if $g=g'$ and $hh'\in E(H)$, or $h=h'$ and $gg'\in E(G)$. Clearly, the Cartesian product is commutative, that is, $G\square H$ is isomorphic to $H\square G$. Moreover, $G\square H$ is $2$-connected whenever $G$ and $H$ are connected. Thus we consider $pvc_k(G\square H)$ for the case $k=2$ in this section. Remind that $d_G(u,v)$ is the shortest distance between the two vertices $u$ and $v$ in graph $G$.

\begin{lem}\label{lemC1}\cite{HIK} Let $(g,h)$ and $(g',h')$ be two vertices of $G\square H$. Then $$d_{G\square H}((g,h),(g',h'))=d_G(g,g')+d_H(h,h').$$
\end{lem}

\begin{thm}\label{thmC1} Let $G$ and $H$ be two nontrivial connected graphs.

$(1)$ If both $G$ and $H$ are complete, then $pvc(G\square H)=1$; otherwise, we have $pvc(G\square H)=2$.

$(2)$ If $G$ and $H$ are two complete graphs of order at least 3, then $pvc_2(G\square H)=1$; otherwise, we have $pvc_2(G\square H)=2$.

$(3)$ If both $G$ and $H$ are complete, then $spvc(G\square H)=1$; otherwise, we have $spvc(G\square H)\\ \leq \min\{spvc(G)\times\chi(H), spvc(H)\times\chi(G)\}$.
\end{thm}

\pf (1) If both $G$ and $H$ are complete, then $diam(G\square H)=2$ and so $pvc(G\square H)=1$ by Proposition \ref{pro1}. Otherwise, we have $diam(G\square H)\geq3$ and so $pvc(G\square H)=2$ by Theorem \ref{thm1}.

(2) First suppose that $G$ and $H$ are two complete graphs of order at least 3. Then $diam(G\square H)=2$ and so $pvc_2(G\square H)\geq1$. Color all the vertices of $G\square H$ with color 1. Next we just need to show that for any two vertices $(g,h)$ and $(g',h')$ in $G\square H$, there exist two vertex-proper paths between them. If $g=g'$, then $(g,h)(g',h')$ and $(g,h)(g,h_0)(g',h')$ are the desired paths, where $h_0\in V(H)\backslash\{h,h'\}$. Similarly, we can get the case that $h=h'$. Now we may assume that $g\neq g'$ and $h\neq h'$. Then $(g,h)(g,h')(g',h')$ and $(g,h)(g',h)(g',h')$ are the desired paths. Hence $pvc_2(G\square H)\leq1$ and so $pvc_2(G\square H)=1$.

Afterwards suppose that $G=K_2$ and $H=K_n$, where $V(G)=\{g_1,g_2\}$. For two vertices $(g_1,h)$ and $(g_2,h)$ of $G\square H$, the edge $(g_1,h)(g_2,h)$ is one desired path but the length of the other desired path is at least 3. Thus $pvc_2(G\square H)\geq2$ and so it remains to show that $pvc_2(G\square H)\leq2$. Define a 2-coloring of $G\square H$ by coloring the vertex $(g_i,h)$ with color $i$ where $i\in \{1,2\}$ and $(g_i,h)\in V(G\square H)$. It is easy to check that there exist two vertex-proper paths between any two vertices in $G\square H$ and so $pvc_2(G\square H)\leq2$. Thus $pvc_2(G\square H)=2$.

Finally we may assume that $G$ is not complete without loss of generality. Then we have $diam(G\square H)\geq3$ and so $pvc_2(G\square H)\geq2$. Next we just need to show that $pvc_2(G\square H)\leq2$. Let $S$ and $T$ be spanning trees of $G$ and $H$, respectively. Then $S\square T$ is a spanning subgraph of $G\square H$ and so it suffices to show that $pvc_2(S\square T)\leq2$. Let $g_0$ and $h_0$ be the roots of $S$ and $T$, respectively. Define a 2-coloring of the vertices of $S\square T$ as follows: For each vertex $(g,h)\in V(S\square T)$, if $d_S(g,g_0)$ and $d_T(h,h_0)$ are of the same parity, then color the vertex $(g,h)$ with color 1; otherwise, color it with color 2. Now it remains to check that there are two vertex-proper paths between any two vertices $(g,h),(g',h')$ in $S\square T$. For any two vertices in $S$ or $T$, there is a path connecting them. Let $gg_1...g_kg'$ and $hh_1...h_lh'$ be two paths from $g$ to $g'$ in $S$ and from $h$ to $h'$ in $T$, respectively. If $g=g'$ and $gg^*$ is an edge in $S$, then $(g,h)(g,h_1)...(g,h_l)(g',h')$ and $(g,h)(g^*,h)(g^*,h_1)...(g^*,h_l)(g^*,h')(g',h')$ are the desired paths. The same is true for the case that $h=h'$. Now we may assume that $g\neq g'$ and $h\neq h'$. Then $(g,h)(g,h_1)...(g,h_l)(g,h')(g_1,h')...(g_k,h')(g',h')$ and $(g,h)(g_1,h)...(g_k,h)(g',h)(g',h_1)...(g',h_l)(g',h')$ are the desired paths. Thus $pvc_2(G\square H)=2$.

(3) If both $G$ and $H$ are complete, then $diam(G\square H)=2$ and so $spvc(G\square H)=1$ by Proposition \ref{pro2}. Otherwise, we will show that $spvc(G\square H)\leq spvc(G)\times\chi(H)$. Firstly define a vertex-coloring $c$ of $G\square H$ with $spvc(G)\times\chi(H)$ colors as follows. We give $G$ a vertex-coloring $c_G$ using $\{1,2,...,spvc(G)\}$ such that $G$ is strong proper vertex-connected, and give $H$ a proper coloring $c_H$ using $\{1,2,...,\chi(H)\}$. For $(g,h)\in V(G\square H)$, where $g\in V(G)$ and $h\in V(H)$, we set $c(g,h)=(c_G(g),c_H(h))$. By this way, we get a vertex-coloring of $G\square H$ with $spvc(G)\times\chi(H)$ colors and it remains to check that for any two vertices $(g,h),(g',h')$ of $G\square H$, there exists a vertex-proper geodesic between them. Let $P=gg_1...g_kg'$ be a vertex-proper geodesic from $g$ to $g'$ in $G$ and $Q=hh_1...h_lh'$ be a shortest path from $h$ to $h'$ in $H$. By Lemma \ref{lemC1}, the path $(g,h)(g_1,h)...(g_k,h)(g_k,h_1)(g_k,h_2)...(g_k,h_l)\\(g_k,h')(g',h')$ is the desired geodesic. Thus, $spvc(G\square H)\leq spvc(G)\times\chi(H)$. By the commutativity of the Cartesian product, we can also deduce that $spvc(G\square H)\leq spvc(H)\times\chi(G)$. Therefore, $spvc(G\square H)\leq \min\{spvc(G)\times\chi(H), spvc(H)\times\chi(G)\}$. \qed

\section{The lexicographic product}

The {\it lexicographic product} $G\circ H$ of graphs $G$ and $H$ is the graph with vertex set $V(G)\times V(H)$, in which two vertices $(g,h),(g',h')$ are adjacent if and only if $gg'\in E(G)$, or $g=g'$ and $hh'\in E(H)$. The lexicographic product is not commutative and is connected whenever $G$ is connected. Moreover, $G\circ H$ is $2$-connected if $G$ and $H$ are connected. Let $d_G(g)$ denote the degree of the vertex $g$ in $G$.

\begin{lem}\label{leml1}\cite{HIK} Let $(g,h)$ and $(g',h')$ be two vertices of $G\circ H$. Then
\begin{eqnarray}d_{G\circ H}((g,h),(g',h'))=
\begin{cases}
d_G(g,g') &if\ g\neq g'; \cr
d_H(h,h') &if\ g=g'\ and\ d_G(g)=0; \cr
\min\{d_H(h,h'),2\} &if\ g=g'\ and\ d_G(g)\neq0.
\end{cases}
\end{eqnarray}
\end{lem}

\begin{thm}\label{thml1} Let $G$ be a nontrivial connected graph and $H$ be a nontrivial graph.

$(1)$ If both $G$ and $H$ are complete, then $pvc(G\circ H)=0$; if $diam(G)\geq3$, then $pvc(G\circ H)=2$; otherwise, we have $pvc(G\circ H)=1$.

$(2)$ If both $G$ and $H$ are complete, then $spvc(G\circ H)=0$; if $diam(G)\geq3$, then $spvc(G\circ H)=2$; otherwise, we have $spvc(G\circ H)=1$.

$(3)$ Let $H$ be a connected graph. If $diam(G)\geq3$, then $pvc_2(G\circ H)=2$; otherwise, we have $pvc_2(G\circ H)=1$.
\end{thm}

\pf (1) If both $G$ and $H$ are complete, then $diam(G\circ H)=1$ and so $pvc(G\circ H)=0$ by Proposition \ref{pro1}. If $G$ is complete and $H$ is not complete, then $diam(G\circ H)=2$ by Lemma \ref{leml1} and so $pvc(G\circ H)=1$ by Proposition \ref{pro1}. Now we may assume that $G$ is not complete. From Lemma \ref{leml1}, it follows that $diam(G\circ H)=diam(G)$. Thus we have that $pvc(G\circ H)=1$ if $diam(G)=2$ and $pvc(G\circ H)=2$ if $diam(G)\geq3$ by Proposition \ref{pro1} and Theorem \ref{thm1}.

(2) If both $G$ and $H$ are complete, then $diam(G\circ H)=1$ and so $spvc(G\circ H)=0$ by Proposition \ref{pro2}. If $G$ is complete and $H$ is not complete, then $diam(G\circ H)=2$ by Lemma \ref{leml1} and so $spvc(G\circ H)=1$ by Proposition \ref{pro2}. Now we may assume that $G$ is not complete. Then $diam(G\circ H)=diam(G)$ by Lemma \ref{leml1}. From Proposition \ref{pro2}, we have that $spvc(G\circ H)=1$ if $diam(G)=2$. Next set $diam(G)\geq3$. Then $spvc(G\circ H)\geq2$ and we just need to show that $spvc(G\circ H)\leq2$. Let $V(H)=\{h_1,h_2,...,h_n\}$. Define a vertex-coloring $c$ of $G\circ H$ with two colors as follows. For $(g,h_i)\in V(G\circ H)$, where $g\in V(G)$ and $i\in[n]$, we set $c(g,h_i)=1$ if $i$ is odd and $c(g,h_i)=2$ if $i$ is even. It suffices to check that there exists a vertex-proper geodesic between any two vertices $(g,h_i),(g',h_j)$ of $G\circ H$. Let $P=gg_1...g_kg'$ be a $g$-$g'$ geodesic in $G$. If $(g,h_i)$ and $(g',h_j)$ are adjacent, then the edge $(g,h_i)(g',h_j)$ is the desired geodesic. Otherwise, if $g=g'$, then set $g^*$ be a neighbor of $g$ in $G$ and so the path $(g,h_i)(g^*,h_1)(g',h_j)$ is the desired geodesic; if $g\neq g'$, then the desired geodesic is $(g,h_i)(g_1,h_1)(g_2,h_2)(g_3,h_1)...(g_k,h_1)(g',h_j)$ when $|P|$ is odd and the desired geodesic is $(g,h_i)(g_1,h_1)(g_2,h_2)(g_3,h_1)...(g_k,h_2)(g',h_j)$ when $|P|$ is even. Thus $spvc(G\circ H)=2$.

(3) If $diam(G)\geq3$, then $diam(G\circ H)=diam(G)$ by Lemma \ref{leml1} and so $pvc_2(G\circ H)\geq2$. Since $G\square H$ is a spanning subgraph of $G\circ H$, $pvc_2(G\circ H)\leq2$ by Theorem \ref{thmC1}$(2)$. Thus $pvc_2(G\circ H)=2$. We now assume that $diam(G)\leq2$. If $diam(G)=diam(H)=1$, then $G\circ H$ is complete and so $pvc_2(G\circ H)=1$ by Lemma \ref{lemJ1}(1). For the other cases, we have $diam(G\circ H)=2$ and so $pvc_2(G\circ H)\geq1$. It suffices to show that $pvc_2(G\circ H)\leq1$. Define a vertex-coloring of $G\circ H$ by coloring each vertex with color 1. Next it remains to check that there are two vertex-proper paths between any two vertices $(g,h),(g',h')$ in $G\circ H$. If $g=g'$, then the paths $(g,h)(g^*,h)(g',h')$ and $(g,h)(g^*,h')(g',h')$ are the desired paths, where $g^*$ is a neighbor of $g$ in $G$. If $h=h'$ and $gg'\in E(G)$, then the edge $(g,h)(g',h')$ and the path $(g,h)(g,h^*)(g',h')$ are the desired paths,where $h^*$ is a neighbor of $h$ in $H$. If $h=h'$ and $gg'\notin E(G)$, then $g$ and $g'$ must have a common neighbor, say $g^*$, since $diam(G)\leq2$. Then the paths $(g,h)(g^*,h)(g',h')$ and $(g,h)(g^*,h^*)(g',h')$ are the desired paths where $h^*\in V(H)\backslash\{h\}$. Now we may assume that $g\neq g'$ and $h\neq h'$. If $gg'\in E(G)$, then the edge $(g,h)(g',h')$ and the path $(g,h)(g,h^*)(g',h')$ are the desired paths, where $h^*$ is a neighbor of $h$ in $H$. Otherwise we have $gg'\notin E(G)$ and then $g$ and $g'$ have a common neighbor, say $g^*$, since $diam(G)\leq2$. Thus the paths $(g,h)(g^*,h)(g',h')$ and $(g,h)(g^*,h')(g',h')$ are the desired paths. Hence $pvc_2(G\circ H)=1$.
\qed

\section{The strong product}

The {\it strong product} $G\boxtimes H$ of graphs $G$ and $H$ is the graph with vertex set $V(G)\times V(H)$, in which two vertices $(g,h),(g',h')$ are adjacent whenever $gg'\in E(G)$ and $h=h'$, or $g=g'$ and $hh'\in E(H)$, or $gg'\in E(G)$ and $hh'\in E(H)$. If an edge of $G\boxtimes H$ belongs to one of the first two types, then we call such an edge a {\it Cartesian edge} and an edge of the last type is called a {\it noncartesian edge}. (The name is due to the fact that if we consider only the first two types, we get the Cartesian product of graphs.) The strong product is commutative and is $2$-connected if $G$ and $H$ are connected. The vertex-connectivity of the strong product was solved recently by Spacapan \cite{S}.

\begin{lem}\label{lemS1}\cite{HIK} Let $(g,h)$ and $(g',h')$ be two vertices of $G\boxtimes H$. Then $$d_{G\boxtimes H}((g,h),(g',h'))=\max \{d_G(g,g'),d_H(h,h')\}.$$
\end{lem}

\begin{thm}\label{thmS1} Let $G$ and $H$ be two nontrivial connected graphs. Then

$(1)$ If both $G$ and $H$ are complete, then $pvc(G\boxtimes H)=0$; if $diam(G)\geq3$ or $diam(H)\geq3$, then $pvc(G\boxtimes H)=2$; otherwise, we have $pvc(G\boxtimes H)=1$.

$(2)$ If $diam(G)\leq2$ and $diam(H)\leq2$ except for the case that $diam(G)=diam(H)=2$ and there exist two vertices with distance $2$ in $G$ and in $H$ have only one common neighbor respectively, then $pvc_2(G\boxtimes H)=1$; otherwise, we have $pvc_2(G\boxtimes H)=2$.

$(3)$ If both $G$ and $H$ are complete, then $spvc(G\boxtimes H)=0$; if either $diam(G)=2$ and $diam(H)\leq2$ or $diam(G)\leq2$ and $diam(H)=2$, then $spvc(G\boxtimes H)=1$; otherwise if $diam(G)\leq2$ and $diam(H)\geq3$, then $spvc(G\boxtimes H)\leq spvc(H)$ and the bound is sharp, and if $diam(G)\geq3$ and $diam(H)\geq3$, then we have $spvc(G\boxtimes H)\leq spvc(G)\times spvc(H)$.
\end{thm}

\pf (1) If both $G$ and $H$ are complete, then $diam(G\boxtimes H)=1$ and so $pvc(G\boxtimes H)=0$ by Proposition \ref{pro1}. If $diam(G)\geq3$ or $diam(H)\geq3$, then $diam(G\boxtimes H)\geq3$ by Lemma \ref{lemS1} and so $pvc(G\boxtimes H)=2$ by Theorem \ref{thm1}. For the other cases, we can deduce that $diam(G\boxtimes H)=2$ and so $pvc(G\boxtimes H)=1$ by Proposition \ref{pro1}.

(2) First suppose that $diam(G)=1$ and $diam(H)\leq2$ by symmetry. If $diam(G)=diam(H)=1$, then $diam(G\boxtimes H)=1$ by Lemma \ref{lemS1} and so $pvc_2(G\boxtimes H)=1$ by Lemma \ref{lemJ1}. Otherwise we have $diam(G)=1$ and $diam(H)=2$. Then $diam(G\boxtimes H)=2$ by Lemma \ref{lemS1} and so $pvc_2(G\boxtimes H)\geq1$. Now color each vertex of $G\boxtimes H$ with color 1. We can always find two disjoint paths of length at most 2 between any two vertices in $G\boxtimes H$. Thus any two vertices of $G\boxtimes H$ are connected by two disjoint vertex-proper paths and then $pvc_2(G\boxtimes H)\leq1$. Thus $pvc_2(G\boxtimes H)=1$. In a similar way, we can deduce $pvc_2(G\boxtimes H)=1$ for the case that $diam(G)=diam(H)=2$ and any two vertices with distance $2$ in $G$ (or $H$) have at least two common neighbors.

Next we consider the case that $diam(G)=diam(H)=2$ and there exist two vertices with distance $2$ in $G$ and in $H$ have only one common neighbor respectively. Let $g$ and $g'$ be the two vertices in $G$ and $g_0$ be their unique common neighbor. Also let $h$ and $h'$ be the two vertices in $H$ and $h_0$ be their unique common neighbor. We can see that the path $(g,h)(g_0,h_0)(g',h')$ is the only path of length 2 between $(g,h)$ and $(g',h')$. Then $pvc_2(G\boxtimes H)\geq2$. Since $G\square H$ is a spanning subgraph of $G\boxtimes H$, we have $pvc_2(G\boxtimes H)\leq pvc_2(G\square H)=2$ by Theorem \ref{thmC1}(2). Thus $pvc_2(G\boxtimes H)=2$.

Finally, suppose that $diam(G)\geq3$ without loss of generality. From Lemma \ref{lemS1}, it follows that $diam(G\boxtimes H)\geq diam(G)\geq3$ and so $pvc_2(G\boxtimes H)\geq2$. Similar to the case above, we have $pvc_2(G\boxtimes H)\leq pvc_2(G\square H)=2$ by Theorem \ref{thmC1}(2) since $G\square H$ is a spanning subgraph of $G\boxtimes H$. Thus $pvc_2(G\boxtimes H)=2$.

(3) If both $G$ and $H$ are complete, then $diam(G\boxtimes H)=1$ by Lemma \ref{lemS1} and so $spvc(G\boxtimes H)=0$ by Proposition \ref{pro2}. If either $diam(G)=2$ and $diam(H)\leq2$ or $diam(G)\leq2$ and $diam(H)=2$, then $diam(G\boxtimes H)=2$ by Lemma \ref{lemS1} and so $spvc(G\boxtimes H)=1$ by Proposition \ref{pro2}.

Afterwards assume that $diam(G)\leq2$ and $diam(H)\geq3$. Next we show that $spvc(G\boxtimes H)\leq spvc(H)$. Firstly define a vertex-coloring $c$ of $G\boxtimes H$ with $spvc(H)$ colors as follows. We give $H$ a vertex-coloring $c_H$ using $\{1,2,...,spvc(H)\}$ such that $H$ is strong proper vertex-connected. For $(g,h)\in V(G\boxtimes H)$, where $g\in V(G)$ and $h\in V(H)$, we set $c(g,h)=c_H(h)$. Then it suffices to check that there exists a vertex-proper geodesic between any two vertices $(g,h),(g',h')$ of $G\boxtimes H$. Let $P=hh_1...h_kh'$ be a vertex-proper $h$-$h'$ geodesic in $H$. Since $diam(G)\leq2$, if $gg'\notin E(G)$, then $g$ and $g'$ must have a common neighbor, say $g^*$. If $g=g'$, then the path $(g,h)(g,h_1)...(g,h_k)(g',h')$ is the desired geodesic. For the case that $h=h'$, if $gg'\in E(G)$, then the edge $(g,h)(g',h')$ is the desired geodesic; otherwise the path $(g,h)(g^*,h)(g',h')$ is the desired geodesic. Now we may assume that $g\neq g'$ and $h\neq h'$. If $gg'\in E(G)$, then the path $(g,h)(g,h_1)...(g,h_k)(g',h')$ is the desired geodesic; otherwise the path $(g,h)(g^*,h_1)(g^*,h_2)...(g^*,h_k)(g',h')$ is the desired geodesic. Thus we have $spvc(G\boxtimes H)\leq spvc(H)$. To show the sharpness of the bound, we consider the example that $G$ is a graph with $diam(G)=2$ and $H$ is a path of length at least 3. By Lemma \ref{lemS1}, it follows that $diam(G\boxtimes H)\geq3$ and so $spvc(G\boxtimes H)\geq2$. Moreover, $spvc(G\boxtimes H)\leq spvc(H)=2$. Hence $spvc(G\boxtimes H)=spvc(H)$.

Finally we consider the case that $diam(G)\geq3$ and $diam(H)\geq3$. To show that $spvc(G\boxtimes H)\leq spvc(G)\times spvc(H)$, we provide a vertex-coloring $c$ of $G\boxtimes H$ with $spvc(G)\times spvc(H)$ colors such that $G\boxtimes H$ is strong proper-vertex connected. First we give $G$ and $H$ two vertex-colorings $c_G$ and $c_H$ using $\{1,2,...,spvc(G)\}$ and $\{1,2,...,spvc(H)\}$ such that $G$ and $H$ are strong proper vertex-connected, respectively. For $(g,h)\in V(G\boxtimes H)$, where $g\in V(G)$ and $h\in V(H)$, we set $c(g,h)=(c_G(g),c_H(h))$. By this way, we get a vertex-coloring $c$ of $G\boxtimes H$ with $spvc(G)\times spvc(H)$ colors and it suffices to check that there exists a vertex-proper geodesic between any two vertices $(g,h),(g',h')$ of $G\boxtimes H$. Let $P=gg_1...g_kg'$ and $Q=hh_1...h_lh'$ be two vertex-proper geodesics from $g$ to $g'$ in $G$ and from $h$ to $h'$ in $H$, respectively. For $h=h'$, if $gg'\in E(G)$, then the edge $(g,h)(g',h')$ is the desired geodesic; otherwise the path $(g,h)(g_1,h)...(g_k,h)(g',h')$ is the desired geodesic. For $h\neq h'$, if $g=g'$ or $gg'\in E(G)$, then the path $(g,h)(g,h_1)...(g,h_l)(g',h')$ is the desired geodesic; otherwise, let $k\leq l$ without loss of generality and then the path $(g,h)(g_1,h_1)...(g_k,h_k)(g_k,h_{k+1})...(g_k,h_l)(g',h')$ is the desired geodesic. Thus we get that $spvc(G\boxtimes H)\leq spvc(G)\times spvc(H)$.\qed

\section{The direct product}

The {\it direct product} $G\times H$ of graphs $G$ and $H$ is the graph with vertex set $V(G)\times V(H)$, in which two vertices $(g,h),(g',h')$ are adjacent if the projections on both coordinates are adjacent, i.e., $gg'\in E(G)$ and $hh'\in E(H)$. It is clearly commutative and associativity also follows quickly. For more general properties we recommend \cite{HIK}. The direct product is the most natural graph product in the sense of categories. But this also seems to be the reason that it is, in general, also the most elusive product of all standard products. For example, $G\times H$ need not be connected even when both factors are. To gain connectedness of $G\times H$ at least one factor must additionally be nonbipartite as shown by Weichsel \cite{W}. Also, the distance formula below for the direct product is far more complicated as it is for other standard products. Here $d_{G}^{e}(g,g')$ represents the length of a shortest even walk between $g$ and $g'$ in $G$, and $d_{G}^{o}(g,g')$ the length of a shortest odd walk between $g$ and $g'$ in $G$. There is no final solution for the connectivity of the direct product, only some partial results are known (see \cite{BS,GV}). But the edge-connectivity of direct products is completely solved in \cite{S1}.

\begin{lem}\label{lemD1}\cite{GH,K} Let $(g,h)$ and $(g',h')$ be two vertices of $G\times H$. Then $$d_{G\times H}((g,h),(g',h'))=\min\{\max \{d_G^{e}(g,g'),d_H^{e}(h,h')\},\max\{d_G^{o}(g,g'),d_H^{o}(h,h')\}\}.$$
\end{lem}

\begin{lem}\label{lemD2} Let $G$ and $H$ be two nontrivial connected graphs. If $G$ or $H$ is nonbipartite, then $diam(G\times H)=2$ if and only if $diam(G)\leq2,diam(H)\leq2$ and each edge of $G$ and $H$ is contained in a triangle in $G$ and $H$ respectively.
\end{lem}

\pf Firstly we prove its sufficiency. Let $(g,h)$ and $(g',h')$ be two vertices in $G\times H$. Since $diam(G)\leq2,diam(H)\leq2$ and each edge of $G$ and $H$ is contained in a triangle in $G$ and $H$ respectively, if $g=g'$, then $d_G^e(g,g')=0$ and $d_H^e(h,h')=2$; if $h=h'$, then $d_G^e(g,g')=2$ and $d_H^e(h,h')=0$; if $g\neq g'$ and $h\neq h'$, then $d_G^e(g,g')=d_H^e(h,h')=2$. Thus $\max \{d_G^{e}(g,g'),d_H^{e}(h,h')\}=2$.  From Lemma \ref{lemD1}, it follows that $d_{G\times H}((g,h),(g',h'))\leq2$ and so $diam(G\times H)\leq2$. For $g=g'$, we have $d_G^o(g,g')\geq3$ and then $d_{G\times H}((g,h),(g',h'))=2$ by Lemma \ref{lemD1}. Hence $diam(G\times H)=2$.

Now we show its necessity. By the definition of direct product, we have $diam(G\times H)\geq\max\{diam(G),diam(H)\}$. Moreover, $diam(G\times H)=2$. Thus $diam(G)\leq2$ and $diam(H)\leq2$. For $(g,h),(g',h)\in V(G\times H)$, where $gg'\in E(G)$, we have $d_H^o(h,h)\geq3$ and $d_{G\times H}((g,h),(g',h))\leq diam(G\times H)=2$. By Lemma \ref{lemD1}, we have $d_G^e(g,g')\leq2$ and then $d_G^e(g,g')=2$ since $gg'\in E(G)$. Thus the edge $gg'$ is contained in a triangle in $G$. Similarly, we can deduce that each edge of $H$ is contained in a triangle in $H$.\qed

For any two vertices $(g,h),(g',h')\in V(G\times H)$, if $g=g'$ or $h=h'$, then there does not exist an edge connecting them. Thus $diam(G\times H)\geq 2$. Moreover, from Propositions \ref{pro1} and \ref{pro2}, Theorem \ref{thm1} and Lemma \ref{lemD2}, we have the following result.

\begin{thm}\label{thmD1} Let $G$ and $H$ be two nontrivial connected graphs, and at least one of them is nonbipartite. Then $pvc(G\times H)=spvc(G\times H)=1$ if $diam(G)\leq2,diam(H)\leq2$ and each edge of $G$ and $H$ is contained in a triangle in $G$ and $H$ respectively, and $pvc(G\times H)=2$ otherwise.
\end{thm}

\begin{thm}\label{thmD2} Let $n$ and $m$ be two integers with $n\geq3$ and $m\geq3$. Then,

$(1)\ spvc(K_n\times K_m)=pvc_2(K_n\times K_m)=1$.

$(2)\ spvc(K_2\times K_m)=pvc_2(K_2\times K_m)=2$.
\end{thm}

\pf $(1)$ It is easy to see that $diam(K_n\times K_m)=2$. By Proposition \ref{pro2}, it follows that $spvc(K_n\times K_m)=1$. Moreover, $pvc_2(K_n\times K_m)\geq1$. It suffices to show that $pvc_2(K_n\times K_m)\leq1$. Define a vertex-coloring of $G\times H$ by coloring each vertex with color 1. Next it remains to check that there are two vertex-proper paths between any two vertices $(g,h),(g',h')$ in $G\times H$. If $g=g'$, then the paths $(g,h)(g_1,h^*)(g',h')$ and $(g,h)(g_2,h^*)(g',h')$ are the desired paths, where $g_1,g_2\in V(G)\backslash\{g\}$ and $h^*\in V(H)\backslash\{h,h'\}$. The same is true for the case that $h=h'$. Now we may assume that $g\neq g'$ and $h\neq h'$. Then the edge $(g,h)(g',h')$ and the path $(g,h)(g^*,h^*)(g',h')$ are the desired paths, where $g^*\in V(G)\backslash\{g,g'\}$ and $h^*\in V(H)\backslash\{h,h'\}$. Hence $pvc_2(K_n\times K_m)=1$.

$(2)$ Since $diam(K_2\times K_m)=3$, $spvc(K_2\times K_m)\geq2$ and $pvc_2(K_2\times K_m)\geq2$. From \cite[Proposition 4]{MYWY}, we have $srvc(K_2\times K_m)=2$. Recall that $spvc(K_2\times K_m)\leq srvc(K_2\times K_m)$ and then $spvc(K_2\times K_m)\leq2$. Thus $spvc(K_2\times K_m)=2$. Next we just need to show that $pvc_2(K_2\times K_m)\leq2$. Define a 2-coloring of $K_2\times K_m$ by coloring the vertex $(g_i,h)$ with color $i$ where $i\in \{1,2\}$ and $(g_i,h)\in V(K_2\times K_m)$. Now it remains to check that there are two vertex-proper paths between any two vertices $(g,h),(g',h')$ in $K_2\times K_m$. If $g=g'$, then the paths $(g,h)(g^*,h^*)(g',h')$ and $(g,h)(g^*,h')(g,h^*)(g^*,h)(g',h')$ are the desired paths, where $g^*\in V(K_2)\backslash\{g\}$ and $h^*\in V(K_m)\backslash\{h,h'\}$. If $h=h'$, then the paths $(g,h)(g',h_1)(g,h_2)(g',h')$ and $(g,h)(g',h_2)(g,h_1)(g',h')$ are the desired paths, where $h_1,h_2\in V(K_m)\backslash\{h\}$. Now we may assume that $g\neq g'$ and $h\neq h'$. Then the edge $(g,h)(g',h')$ and the path $(g,h)(g',h^*)(g,h')(g',h)(g,h^*)(g',h')$ are the desired paths, where $h^*\in V(K_m)\backslash\{h,h'\}$. Hence $pvc_2(K_2\times K_m)=2$.\qed

\begin{thm}\label{thmD3} Let $H$ be a nontrivial connected graph with $diam(H)\geq2$. Then we have

$(1)$ If $H$ is nonbipartite, then $spvc(K_2\times H)=2$.

$(2)$ For $n\geq4$, if $diam (H)=2$ and each edge of $H$ is contained in a triangle in $H$, then $spvc(K_n\times H)=1$; otherwise, we have $spvc(K_n\times H)=2$.

$(3)$ If $diam (H)=2$ and each edge of $H$ is contained in a triangle in $H$, then $spvc(K_3\times H)=1$; otherwise, we have $spvc(K_3\times H)=2$ when $H$ is a tree and $2\leq spvc(K_3\times H)\leq3$ when $H$ is not.
\end{thm}

\pf $(1)$ Since $diam(K_2\times H)\geq3$, $spvc(K_2\times H)\geq2$ and so we just need to show that $spvc(K_2\times H)\leq2$. Define a 2-coloring of $K_2\times H$ by coloring the vertex $(g_i,h)$ with color $i$ where $i\in \{1,2\}$ and $(g_i,h)\in K_2\times H$. Now it remains to check that there is a vertex-proper geodesic between any two vertices $(g,h),(g',h')$ in $K_2\times H$. Let $hh_1...h_sh'$ be a shortest odd walk and $hh'_1...h'_th'$ be a shortest even walk between $h$ and $h'$ in $H$. If $g=g'$, then we can get that $d_{G\times H}((g,h),(g',h'))=d_H^{e}(h,h')\geq2$ by Lemma \ref{lemD1} and so the path $(g,h)(g^*,h'_1)...(g^*,h'_t)(g',h')$ is the desired geodesic, where $g^*\in V(K_2)\backslash\{g\}$. If $h=h'$, then we can get that $d_{G\times H}((g,h),(g',h'))=d_H^{o}(h,h')\geq3$ by Lemma \ref{lemD1} and so the path $(g,h)(g',h_1)...(g,h_s)(g',h')$ is the desired geodesic. Now we may suppose that $g\neq g'$ and $h\neq h'$. From Lemma \ref{lemD1}, we can deduce that $d_{G\times H}((g,h),(g',h'))=d_H^{o}(h,h')$. If $hh'\in E(H)$, then the edge $(g,h)(g',h')$ is the desired geodesic; otherwise, the path $(g,h)(g',h_1)...(g,h_s)(g',h')$ is the desired geodesic. Hence $spvc(K_2\times H)=2$.

$(2)$ If $diam (H)=2$ and each edge of $H$ is contained in a triangle in $H$, then $diam(K_n\times H)=2$ by Lemma \ref{lemD2} and so $spvc(K_n\times H)=1$ by Proposition \ref{pro2}. Otherwise, we have $diam(K_n\times H)\geq3$ and so $spvc(K_n\times H)\geq2$. Next we just need to show that $spvc(K_n\times H)\leq2$. Let $V(K_n)=\{g_1,g_2,...,g_n\}$. Define a vertex-coloring $c$ of $K_n\times H$ with two colors as follows. For $(g_i,h)\in V(K_n\times H)$, where $i\in[n]$ and $h\in V(H)$, we set $c(g_i,h)=1$ if $i$ is odd and $c(g_i,h)=2$ if $i$ is even. It suffices to check that there exists a vertex-proper geodesic between any two vertices $(g_i,h),(g_j,h')$ of $K_n\times H$. Let $hh_1...h_sh'$ be a shortest odd walk and $hh'_1...h'_th'$ be a shortest even walk between $h$ and $h'$ in $H$.

First suppose that $i=j$. If $hh'\in E(H)$, then $d_{K_n\times H}((g_i,h),(g_j,h'))=\min\{d_H^{e}(h,h'),\\d_{K_n}^{o}(g_i,g_j)\}$ by Lemma \ref{lemD1}, where $d_H^{e}(h,h')\geq2$ and $d_{K_n}^{o}(g_i,g_j)=3$. Thus if $h$ and $h'$ have a common neighbor, say $h^*$, then $d_{K_n\times H}((g_i,h),(g_j,h'))=d_H^{e}(h,h')=2$ and so the path $(g_i,h)(g_k,h^*)(g_j,h')$ is the desired geodesic, where $g_k\in V(K_n)\backslash\{g_i\}$; otherwise, $d_{K_n\times H}((g_i,h),(g_j,h'))=d_{K_n}^{o}(g_i,g_j)=3$ and the desired geodesic is $(g_i,h)(g_k,h')(g_l,h)(g_j,h')$, where $g_k,g_l\in V(K_n)\backslash\{g_i\}$ and $|k-l|$ is odd. If $hh'\notin E(H)$, then $d_{K_n\times H}((g_i,h),(g_j,h'))=\min\{d_H^{e}(h,h'),d_H^{o}(h,h')\}$ by Lemma \ref{lemD1}, where $d_H^{e}(h,h')\geq2$ and $d_H^{o}(h,h')\geq3$. Thus if $d_H^{e}(h,h')\geq d_H^{o}(h,h')$, then the path $(g_i,h)(g_k,h_1)(g_l,h_2)(g_k,h_3)(g_l,h_4)(g_k,h_5)...(g_l,h_s)(g_j,h')$ is the desired geodesic, where $g_k,g_l\in V(K_n)\backslash\{g_i\}$ and $|k-l|$ is odd; otherwise, the desired geodesic is $(g_i,h)(g_p,h'_1)(g_i,h'_2)...(g_p,h'_t)(g_j,h')$, where $g_p\in V(K_n)\backslash\{g_i\}$ and $|p-i|$ is odd.

Then suppose that $h=h'$. From Lemma \ref{lemD1}, we have $d_{K_n\times H}((g_i,h),(g_j,h'))=d_{K_n}^{e}(g_i,g_j)=2$ and so the desired geodesic is $(g_i,h)(g_k,h^*)(g_j,h')$ where $g_k\in V(K_n)\backslash\{g_i,g_j\}$ and $h^*$ is a neighbor of $h$ in $H$.

Finally we suppose that $i\neq j$ and $h\neq h'$. If $hh'\in E(H)$, then the edge $(g_i,h)(g_j,h')$ is the desired geodesic. Otherwise, we can get that $d_{K_n\times H}((g_i,h),(g_j,h'))=\min\{d_H^{e}(h,h'),\\d_H^{o}(h,h')\}$ by Lemma \ref{lemD1}, where $d_H^{e}(h,h')\geq2$ and $d_H^{o}(h,h')\geq3$. Thus if $d_H^{e}(h,h')\geq d_H^{o}(h,h')$, then the path $(g_i,h)(g_j,h_1)(g_k,h_2)...(g_k,h_s)(g_j,h')$ is the desired geodesic where $g_k\in V(K_n)\backslash\{g_i,g_j\}$ and $|k-j|$ is odd; otherwise, the desired geodesic is $(g_i,h)(g_p,h'_1)(g_j,h'_2)\\...(g_p,h'_t)(g_j,h')$ where $g_p\in V(K_n)\backslash\{g_i,g_j\}$ and $|p-j|$ is odd.

Therefore $spvc(K_n\times H)=2$.

$(3)$ If $diam (H)=2$ and each edge of $H$ is contained in a triangle in $H$, then $diam(K_3\times H)=2$ by Lemma \ref{lemD2} and so $spvc(K_3\times H)=1$ by Proposition \ref{pro2}. Otherwise, we have $diam(K_3\times H)\geq3$ and so $spvc(K_3\times H)\geq2$. Suppose that $H$ is a tree. Then we just need to show that $spvc(K_3\times H)\leq2$. Let $h_0$ be the root of the tree $H$. Define a 2-coloring of the vertices of $K_3\times H$ as follows: For each vertex $(g,h)\in K_3\times H$, if $d_H(h,h_0)$ is odd, then color the vertex $(g,h)$ with color 1; otherwise, color it with color 2. Now it remains to check that there is a vertex-proper geodesic between any two vertices $(g,h),(g',h')$ in $K_3\times H$. For any two vertices in $H$, there is a path connecting them. Let $P=hh_1...h_th'$ be the path from $h$ to $h'$ in $H$. If $h=h'$, then we have $d_{K_3\times H}((g,h),(g',h'))=d_{K_3}^{e}(g,g')=2$ by Lemma \ref{lemD1} and so the path $(g,h)(g^*,h^*)(g',h')$ is the desired geodesic where $g^*\in V(K_3)\backslash\{g,g'\}$ and $h^*$ is a neighbor of $h$ in $H$. Then we consider the case that $g=g'$. If $hh'\in E(H)$, then $d_{K_3\times H}((g,h),(g',h'))=d_{K_3}^{o}(g,g')=3$ by Lemma \ref{lemD1} and so the path $(g,h)(g_1,h')(g_2,h)(g',h')$ is the desired geodesic, where $g_1,g_2\in V(K_3)\backslash\{g\}$. Otherwise, $d_{K_3\times H}((g,h),(g',h'))=\min\{d_H^{e}(h,h'),d_H^{o}(h,h')\}$ by Lemma \ref{lemD1}. Thus if $d_H^{e}(h,h')\geq d_H^{o}(h,h')$, then the path $(g,h)(g_1,h_1)(g_2,h_2)(g,h_3)(g_2,h_4)...(g_2,h_t)(g',h')$ is the desired geodesic where $g_1,g_2\in V(K_3)\backslash\{g\}$; otherwise, the desired geodesic is $(g,h)(g^*,h_1)(g,h_2)\\...(g^*,h_t)(g',h')$ where $g^*\in V(K_3)\backslash\{g\}$. Now we may assume that $g\neq g'$ and $h\neq h'$. If $hh'\in E(H)$, then the edge $(g,h)(g',h')$ is the desired geodesic. Otherwise, we have $d_{K_3\times H}((g,h),(g',h'))=\min\{d_H^{e}(h,h'),d_H^{o}(h,h')\}$ by Lemma \ref{lemD1}. Thus if $d_H^{e}(h,h')\geq d_H^{o}(h,h')$, then the path $(g,h)(g',h_1)(g,h_2)...(g,h_t)(g',h')$ is the desired geodesic; otherwise, the desired geodesic is $(g,h)(g^*,h_1)(g,h_2)...(g^*,h_t)(g',h')$ where $g^*\in V(K_3)\backslash\{g,g'\}$. Hence $spvc(K_3\times H)=2$ when $H$ is a tree.

Next suppose that $H$ is not a tree. Then we just need to show that $spvc(K_3\times H)\leq3$. Let $V(K_3)=\{g_1,g_2,g_3\}$. Define a vertex-coloring $c$ of $K_3\times H$ with three colors as follows. For $(g_i,h)\in V(K_3\times H)$, we set $c(g_i,h)=i$ where $i\in\{1,2,3\}$. Similar to the discussion in $(2)$, we can check that there exists a vertex-proper geodesic between any two vertices $(g_i,h),(g_j,h')$ of $K_3\times H$. Thus $2\leq spvc(K_3\times H)\leq3$.\qed


\begin{thebibliography}{1}

\bibitem{ACKP}
 B. S. Anand, M. Changat, S. Klav$\check{z}$ar, I. Peterin, Convex sets in lexicographic products of graphs, {\it Graphs Combin.} {\bf 28(1)} (2012) 77-84.

\bibitem{ALL}
 E. Andrews, E. Laforge, C. Lumduanhom, P. Zhang, On proper-path colorings
 in graphs, {\it J. Combin. Math. Combin. Comput} {\bf 97} (2016) 189-207.

\bibitem{B}
 J.A. Bondy, U.S.R. Murty, {\it Graph Theory}, GTM 244, Springer, 2008.

\bibitem{BFG}
 V. Borozan, S. Fujita, A. Gerek, C. Magnant, Y. Manoussakis, L. Montero,
 Z. Tuza, Proper connection of graphs, {\it Discrete Math.} {\bf 312(17)} (2012) 2550-2560.

\bibitem{BS}
 B. Bre$\check{s}$ar, S. $\check{S}$pacapan, On the connectivity of the direct product of graphs, {\it Australas. J. Comb.} {\bf 41} (2008) 45-56.

\bibitem{CJM}
 G. Chartrand, G.L. Johns, K.A. McKeon, P. Zhang, Rainbow connection in
 graphs, {\it Math. Bohemica} {\bf 133(1)} (2008) 85-98.

\bibitem{CJMZ}
 G. Chartrand, G.L. Johns, K.A. Mckeon, P. Zhang, The rainbow connectivity
 of a graph, {\it Networks} {\bf 54(2)} (2009) 75-81.

\bibitem{CLS1}
 L. Chen, X. Li, M. Liu, Nordhaus-Gaddum-type theorem for the rainbow
 vertex-connection number of a graph, {\it Util. Math.} {\bf 86} (2011) 335-340.

\bibitem{CLS2}
 L. Chen, X. Li, Y. Shi, The complexity of determining the rainbow
 vertex-connection of a graph, {\it Theoret. Comput. Sci.} {\bf 412(35)} (2011) 4531-4535.

\bibitem{GH}
 A.A. Ghidewon, R. Hammack, Centers of tensor product of graphs, {\it Ars Comb.} {\bf 74} (2005) 201-211.

\bibitem{GMP}
 T. Gologranc, G. Meki$\check{s}$, I. Peterin, Rainbow connection and graph products, {\it Graphs Combin.} {\bf 30(3)} (2014) 591-607.

\bibitem{GV}
 R. Guji, E. Vumar, A note on the connectivity of Kronecker products of graphs, {\it Appl. Math. Lett.} {\bf 22(9)} (2009) 1360-1363.

\bibitem{HIK}
 R. Hammack, W. Imrich, S. Klav$\breve{z}$r, {\it Handbook of product graphs}, 2nd ed., CRC Press, New York, 2011.

\bibitem{JLZZ}
 H. Jiang, X. Li, Y. Zhang, Y. Zhao, On (strong) proper vertex-connection of graphs, {\it Bull. Malays. Math. Sci. Soc.} {\bf 93(1)} (2015) 1-11.

\bibitem{K}
 S.R. Kim, Centers of a tensor composite graph, {\it Congr. Numer.} {\bf 81} (1991) 193-203.

\bibitem{KS}
 S. Klav$\check{z}$ar, S. $\check{S}$pacapan, On the edge-connectivity of Cartesian product graphs, {\it Asian-Eur. J. Math.} {\bf 1} (2008) 93-98.

\bibitem{KY}
M. Krivelevich, R. Yuster, The rainbow connection of a graph is (at
most) reciprocal to its minimum degree, {\it J. Graph Theory} {\bf
63(3)} (2010) 185-191.

\bibitem{LLZ}
 E. Laforge, C. Lumduanhom, P. Zhang, Characterizations of graphs
 having large proper connection numbers, {\it Discuss. Math. Graph Theory} {\bf
36(2)} (2016) 439-453.

\bibitem{LMS1}
 X. Li, Y. Mao, Y. Shi, The strong rainbow vertex-connection of
 graphs, {\it Util. Math.} {\bf 93} (2014) 213-223.

\bibitem{LS}
 X. Li, Y. Shi, On the rainbow vertex-connection, {\it Discuss. Math. Graph Theory} {\bf 33} (2013) 307-313.

\bibitem{LSS}
 X. Li, Y. Shi, Y. Sun, Rainbow connections of graphs: A survey, {\it Graphs Combin.} {\bf 29(1)} (2013) 1-38.

\bibitem{LSu}
 X. Li, Y. Sun, {\it Rainbow Connections of Graphs,} Springer Briefs in Math., Springer, New York, 2012.

\bibitem{LMS}
 H. Liu, $\hat{A}$. Mestre, T. Sousa, Rainbow vertex $k$-connection in graphs, {\it Discrete Appl. Math.} {\bf 162} (2013) 2549-2555.

\bibitem{MYWY}
 Y. Mao, F. Yanling, Z. Wang, C. Ye, Rainbow vertex-connection and graph products, {\it Int. J. Comput. Math.} {\bf 93(7)} (2016) 1078-1092.

\bibitem{MYWY1}
 Y. Mao, F. Yanling, Z. Wang, C. Ye, Proper connection number and graph products, {\it Bull. Malays. Math. Sci. Soc.} DOI 10.1007/s40840-016-0442-z, in press.

\bibitem{NS}
 R. J. Nowakowski, K. Seyffarth, Small cycle double covers of products. I. Lexicographic product with paths and cycles, {\it J. Graph Theory} {\bf 57(2)} (2008) 99-123.

\bibitem{P}
 I. Peterin, Intervals and convex sets in strong product of graphs, {\it Graphs Combin.} {\bf 29(3)} (2013) 705-714.

\bibitem{S}
 S. $\check{S}$pacapan, Connectivity of strong products of graphs, {\it Graphs Combin.} {\bf 26(3)} (2010) 457-467.

\bibitem{S1}
 S. $\check{S}$pacapan, A characterization of edge-connectivity of direct products of graphs, {\it Discret. Math.} {\bf 313(12)} (2013) 1385-1393.

\bibitem{W}
 P.M. Weichsel, The Kronecker product of graphs, {\it Proc. Am. Math. Soc.} {\bf 13(1)} (1962) 47-52.

\bibitem{Z}
 X. Zhu, Game coloring the Cartesian product of graphs, {\it J. Graph Theory} {\bf 59(4)} (2008) 261-278.

\end{thebibliography}
\end{document}